\pdfoutput=1

\documentclass[letterpaper, 10 pt, conference]{ieeeconf}  

\IEEEoverridecommandlockouts                              

\usepackage{amsmath,amssymb,amsfonts}
\usepackage{mathtools}
\mathtoolsset{showonlyrefs}
\usepackage{hyperref}
\usepackage{url}
\usepackage{graphicx}
\usepackage{xcolor}
\usepackage{subcaption}
\usepackage{multirow}
\usepackage{array}
\usepackage{cite}
\newcolumntype{C}[1]{>{\centering\let\newline\\\arraybackslash\hspace{0pt}}m{#1}}
\newtheorem{thm}{Theorem}
\newtheorem{lem}[thm]{Lemma}

\newtheorem{cor}[thm]{Corollary}
\newtheorem{defn}{Definition}

\newtheorem{rem}{Remark}
\newtheorem{ass}{Assumption}

\usepackage{ifthen}
\newboolean{online}
\setboolean{online}{true}

\title{\LARGE \bf
Dissipative Gradient Descent Ascent Method: A Control Theory Inspired Algorithm for Min-max Optimization
}

\author{Tianqi Zheng$^{1}$, Nicolas Loizou$^{1}$, Pengcheng You$^{2}$ and Enrique Mallada$^{1}$
\thanks{$^{1}$T. Zheng, and E. Mallada are with the Department of Electrical and Computer Engineering, N. Loizou is with the Department of Applied Mathematics and Statistics at Johns Hopkins University, Baltimore, MD 21218, USA 
        {\tt\footnotesize \{tzheng8,nloizou,mallada\}@jhu.edu}}%
\thanks{$^{2}$P. You is with the Department of Industrial Engineering and Management at Peking University, Beijing, China 
        {\tt\footnotesize pcyou@pku.edu.cn}}%
}

\begin{document}

\maketitle
\thispagestyle{empty}
\pagestyle{empty}

\begin{abstract}
Gradient Descent Ascent (GDA) methods for min-max optimization problems typically produce oscillatory behavior that can lead to instability, e.g., in bilinear settings.
To address this problem, we introduce a dissipation term into the GDA updates to dampen these oscillations. The proposed Dissipative GDA (DGDA) method can be seen as performing standard GDA on a state-augmented and regularized saddle function that does not strictly introduce additional convexity/concavity. We theoretically show the linear convergence of DGDA in the bilinear and strongly convex-strongly concave settings and assess its performance by comparing DGDA with other methods such as GDA, Extra-Gradient (EG), and Optimistic GDA.
Our findings demonstrate that DGDA surpasses these methods, achieving superior convergence rates. We support our claims with two numerical examples that showcase DGDA's effectiveness in solving saddle point problems.

\end{abstract}

\section{Introduction}
In recent years, there has been a significant focus on solving saddle point problems, namely min-max optimization problems \cite{mokhtari2020unified,azizian2020tight,gorbunov2022stochastic,loizou2021stochastic,beznosikov2023stochastic,gidel2018variational}. These problems have garnered considerable attention, particularly in fields such as Generative Adversarial Networks (GANs) \cite{goodfellow2014generative,gidel2018variational,daskalakis2017training}, Reinforcement Learning (RL) \cite{pfau2016connecting}, and Constrained RL (C-RL) \cite{zheng2022constrained,ding2023last}. However, a major challenge that persists in these approaches is the instability of the training process. That is, solving the min-max optimization problem via running the standard Gradient Descent Ascent (GDA) algorithm often leads to unstable oscillatory behavior rather than convergence to the optimal solution. This is particularly illustrated in bilinear min-max problems, such as the training of Wasserstein GANs \cite{adler2018banach} or solving C-RL problems in the occupancy measure space \cite{altman1999constrained}, for which the standard GDA fails to converge~\cite{mokhtari2020unified,azizian2020tight}.

In order to understand the instability of the GDA method and further tackle its limitation, we draw inspiration from the control-theoretic notions of dissipativity~\cite{sastry2013nonlinear}, which enables the design of stabilizing controllers using dynamic (state-augmented) components that seek to dissipate the energy generated by the unstable process. This aligns with recent work that leverages control theory tools in the analysis and design of optimization algorithms~\cite {nelson2018integral,hu2017dissipativity, fazlyab2018analysis, lessard2016analysis,zhang2021unified}. From a dynamical system point of view, dissipativity theory characterizes the manner in which energy dissipates within the system and drives it towards equilibrium. It provides a direct way to construct a Lyapunov function, which further relates the rate of decrease of this internal energy to the rate of convergence of the algorithm.

We motivate our developments by looking first at a simple scalar bilinear problem wherein, the energy of the system, expressed as the square 2-norm distance to the saddle, strictly increases on every iteration, leading to oscillations of increasing amplitude.  
To tackle this unstable oscillating behavior, we propose the Dissipative GDA method, which, as the name suggests, incorporates a simple friction term to GDA updates to dissipate the internal energy and stabilize the system. 
Our algorithm can be seen as a discrete-time version of \cite{you2021saddle}, which has been applied to solve the C-RL problems \cite{zheng2022constrained}. In this work, we build on this literature, making the following contributions:

\noindent
\emph{1. Novel control theory inspired algorithm:}  
We illustrate how to use control theoretic concepts of dissipativity theory to design an algorithm that can stabilize the unstable behavior of GDA. Particularly, we show that by introducing a friction term, the proposed DGDA algorithm dissipates the stored internal energy and converges toward equilibrium.

\noindent
\emph{2. Theoretical analysis with better rates:} We establish the linear convergence of the DGDA method for bilinear and strongly convex-strongly concave saddle point problems. In both settings, we show that the DGDA method outperforms other state-of-the-art first-order explicit methods, surpassing the standard known linear convergence rate (see Table \ref{table: Comparison of bilinear convergence rate} and \ref{table: Comparison of str convergence rate}).

\noindent
\emph{3. Numerical Validation:} We corroborate our theoretical results with numerical experiments by evaluating the performance of the DGDA method with GDA, EG, and OGDA methods. When applied to solve bilinear and strongly convex-strongly concave saddle point problems, the DGDA method systematically outperforms other methods in terms of convergence rate.

\noindent
\emph{Outline:} The rest of the paper is organized as follows. In Section \ref{Sec: Background}, we provide some preliminary definitions and background. In Section \ref{Sec: Algorithms and Theory}, we leverage tools from dissipativity theory and propose the Dissipative GDA (DGDA) algorithm to tackle the unstable oscillatory behavior of GDA methods. In Section \ref{Sec: Convergence}, we establish its linear convergence rate for bilinear and strongly convex-strongly concave problems, which outperforms state-of-the-art first-order explicit algorithms, including GDA, EG, and OGDA methods. In Section \ref{Sec: Numerical Experiments}, we support our claims with two numerical examples. We close the paper with concluding remarks and future research directions in Section \ref{Sec: Conclusion}.

\begin{table}[t]
\begin{center}
\scalebox{1}{
\begin{tabular}{ C{0.5cm}|C{2cm}C{2cm}C{2cm} }
\hline
Bil. &  Mokhtari, 2020 & Azizian,  2020 & This Work  
\\ \hline 
EG &  $\frac{\kappa^{-1}}{20}$ & $\frac{\kappa^{-1}}{64}$ &  - \\
OG & $\frac{\kappa^{-1}}{800}$ &$\frac{\kappa^{-1}}{128}$ & -  \\
DG  &- &- & $\frac{\kappa^{-1}}{4}$ \\
\end{tabular}}
\end{center}
\caption{Summary of the global convergence results for EG, OGDA, and DGDA methods with bilinear objective functions. If a result shows that the iterates converge as $\mathcal{O}((1 - r)^t)$, the quantity r is reported (the larger the better). $\kappa$ represents the condition number.}
\label{table: Comparison of bilinear convergence rate}
\end{table}

\section{Problem Formulation}\label{Sec: Background}
In this paper, we study the problem of finding saddle points in the min-max optimization problem:
\begin{align}\label{eq: min-max optimization problem}
    \min_{ x \in \mathbb{R}^n}\max_{ y \in \mathbb{R}^m} f( x, y),
\end{align}
where the function $f:\mathbb{R}^n \times \mathbb{R}^m \to \mathbb{R} $ is a convex-concave function. Precisely, $f(\cdot, y)$ is convex for all $ y \in \mathbb{R}^m $ and $f( x,\cdot)$ is concave for all $ x \in \mathbb{R}^n $. We seek to develop a novel optimization algorithm that converges to some saddle point $( x^*,  y^*)$ of Problem \ref{eq: min-max optimization problem}.
\begin{defn}[Saddle Point]
    A point $( x^*,  y^*) \in \mathbb{R}^n \times \mathbb{R}^m$ is a saddle point of convex-concave function \eqref{eq: min-max optimization problem} if and only if it satisfies $f( x^*, y) \leq  f( x^*, y^* ) \leq  f( x, y^* )$    for all $ x \in \mathbb{R}^n ,  y \in \mathbb{R}^m $.
\end{defn}
Throughout this paper, we consider two specific instances of Problem \ref{eq: min-max optimization problem} commonly studied in related literature: strongly convex-strongly concave and bilinear functions. Herein, we briefly present some definitions and properties.
\begin{defn}[Strongly Convex]\label{def: strongly convex concave}
    A differentiable function $f: \mathbb{R}^n \to \mathbb{R}$ is said to be $\mu$-strongly convex if 
    $ f( w) \geq f( w') + \nabla f( w)^T( w- w') + \frac{\mu}{2} \| w -  w' \|^2$. 
\end{defn} 
Notice that if $\mu=0$, then we recover the definition of convexity for a continuously differentiable function and $f( w) $ is $\mu$-strongly concave if $-f( w) $ is $\mu$-strongly convex. Another important property commonly used in the convergence analysis of optimization algorithms is the Lipschitz-ness of the gradient $\nabla f(w)$.

\begin{defn}[$L$-Lipschitz] \label{def: Lipschitz}
A function $F:\mathbb{R}^n \to \mathbb{R}^m$ is L-Lipschitz if $ \forall  w, w'\in\mathbb{R}^n$, we have $\|F( w) - F( w')  \| \leq L \|  w -  w' \|$.
\end{defn}

Combining the above two properties leads to the first important class of problem that has been extensively studied  \cite{rockafellar1976monotone,tseng1995linear,azizian2020tight,mokhtari2020unified}.
\begin{ass}\label{ass: Str} \textit{(Strongly strongly convex-concave functions  with L-Lipschitz Gradient)}
    The function $f:\mathbb{R}^n \times \mathbb{R}^m \to \mathbb{R} $ is continuously differentiable, $\mu$ strongly convex in $ x$, and $\mu$ strongly concave in $ y$. Further, the gradient vector $(\nabla_x f(x,y); -\nabla_y f(x,y))$ is $L$-Lipschitz.
\end{ass}

It is also crucial to consider situations where the objective function is bilinear. Such bilinear min-max problems often appear when solving constrained reinforcement learning problems~\cite{zheng2022constrained,wang2020randomized}, and training of WGANs \cite{adler2018banach}.

\begin{ass}[Bilinear function] \label{ass: Bilinear}
    The function $f:\mathbb{R}^n \times \mathbb{R}^m \to \mathbb{R} $ is a bilinear function if it can be written in the form $f( x, y) =  x ^T   A  y$. For simplicity, we further assume that the matrix $   A \in \mathbb{R}^{m \times n}$ is full rank, with $m\leq n$.
\end{ass}

As seen in Table \ref{table: Comparison of bilinear convergence rate} and \ref{table: Comparison of str convergence rate} as well as in Section \ref{Sec: Convergence}, the linear convergence rates of existing algorithms are frequently characterized by the \textit{condition number} $\kappa$. Specifically, when the objective function is bilinear, the condition number is defined as $\kappa:=\sigma_{\max}^2(A)/\sigma_{\min}^2(A)$, where $\sigma_{\max}( M) $ and $\sigma_{\min}(M)$ denote the largest singular value and smallest singular of a matrix $M$ respectively. When the objective function is strongly convex-strongly concave with the $L$-Lipschitz gradient, the condition number of the problem is defined as $\kappa:=L /\mu$.

\section{Dissipative Gradient Descent Ascent Algorithm}\label{Sec: Algorithms and Theory}
\begin{table}[t]
\begin{center}
\scalebox{0.9}{
\begin{tabular}{ C{0.5cm}|C{1.5cm}C{1.8cm}C{1.8cm}C{1.5cm} }
\hline
S.C & Zhang, 2021 & Mokhtari, 2020 & Azizian, 2020 & This Work  
\\ \hline 
GD & $\kappa^{-2}$ & -& - &- \\
EG & - & $\frac{\kappa^{-1}}{4} $ & $\frac{\kappa^{-1}}{4} +\epsilon$ &  - \\
OG    & - & $\frac{\kappa^{-1}}{4} $ &$\frac{\kappa^{-1}}{4}+\epsilon$ &  - \\
DG   &-& - &- &\!\!\!\! {\small $\kappa^{-1}\!\!\!~-~\mathcal{O}(\kappa^{-2}) $}
\end{tabular}}
\end{center}
\caption{Summary of the global convergence results for GDA, EG, OGDA, and DGDA methods with strongly convex-strongly concave and $L$-Lipschitz objective functions. The table reports the term $r$ of a $(1-r)$ linear rate.
The constant $\epsilon > 0$ depends on the problem.  }
\label{table: Comparison of str convergence rate}
\end{table}

This section introduces the proposed first-order method for solving the min-max optimization problem \ref{eq: min-max optimization problem}. 
The algorithm can be seen as a discretization of the algorithm proposed by \cite{you2021saddle}, where the authors introduce a regularization framework for continuous saddle flow dynamics that guarantees asymptotic convergence to a saddle point under mild assumptions. 
However, the continuous-time analysis presented in \cite{you2021saddle} does generally extend to discrete time. In this work, we will show the linear convergence of the discrete-time version of this algorithm. Moreover, we illustrate a novel control-theoretic design methodology for optimization algorithms that can broadly use dissipativity tools to study the convergence of the resulting algorithm.

Our results build on gaining an intuitive understanding of the problems that one encounters when  applying the vanilla Gradient Descent Ascent method to solve saddle point problems \eqref{eq: min-max optimization problem}:\\
\textbf{Gradient Descent Ascent (GDA)}
\begin{equation}\label{eq: GD}
\begin{aligned}
     &x_{k+1} =  x_k - \eta \nabla_ x f( x_k, y_k), \\
     &y_{k+1} =  y_k + \eta \nabla_ y f( x_k, y_k).
\end{aligned}    
\end{equation}
When \eqref{eq: min-max optimization problem} is strongly convex-strongly concave, and has L-Lipschitz gradients, the GDA method provides linear convergence, with step size $\eta = \mu/L^2$ and a know rate estimate of $ 1 - 1/\kappa^2$ \cite{grimmer2023landscape}. 
However, when the problem is bilinear, the standard GDA method fails to converge.

\subsection{Control Theory-Based Motivation}\label{subsec: Motivation}
\begin{figure}[h]
\centering
\includegraphics[width=\columnwidth]{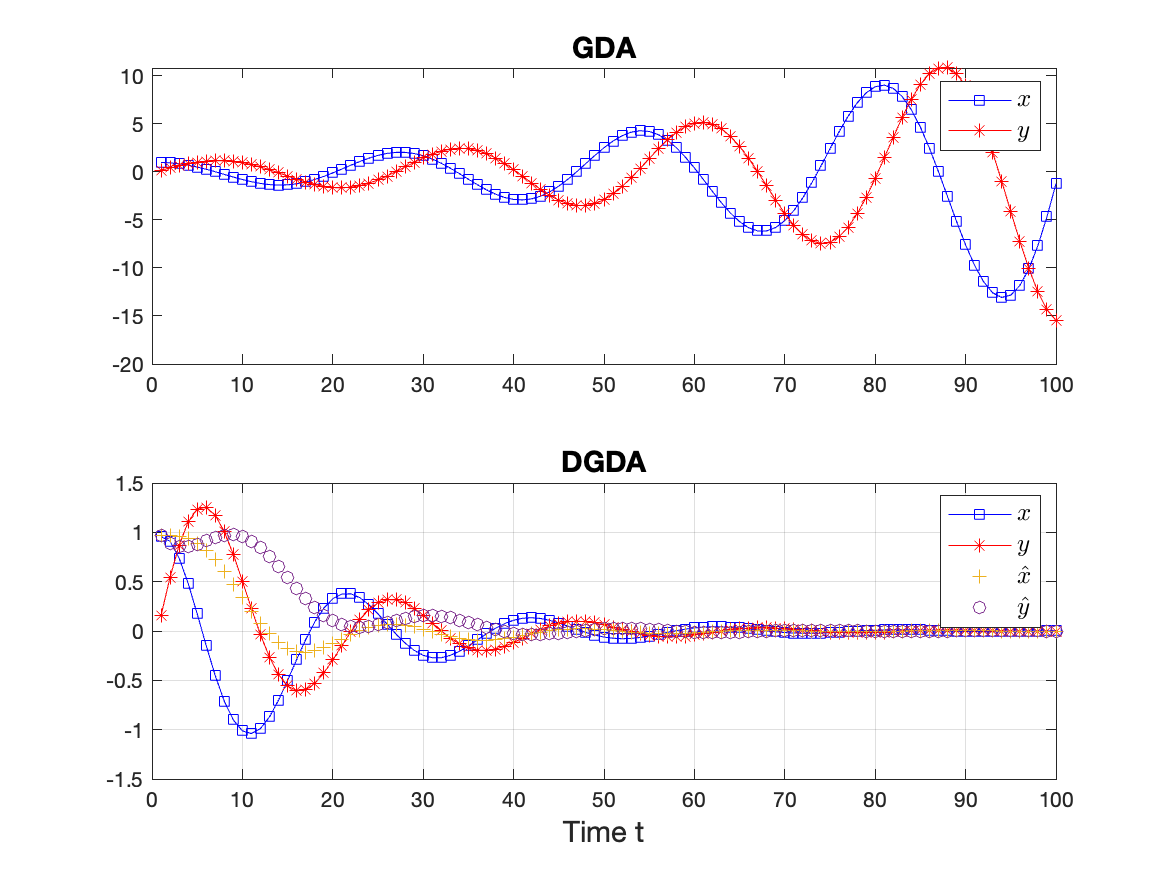}
    \caption{Trajectories of states  for GDA and DGDA for the simple bilinear objective function $f(x,y):=xy$.}
\label{fig:bilinear2d}
\end{figure}

Before we dive into our algorithm, we would like to illustrate our key insight from control theory about how to characterize the energy stored in the system and how to introduce friction to dissipate energy and drive the system toward equilibrium. To get started, consider the following motivating example with a simple bilinear objective function: 
\begin{align}\label{eq: Bilinear Simple}
    \min_{x} \max_{y}f(x,y):=xy
\end{align}
We seek to analyze the behavior of (a possibly controlled) vanilla GDA algorithm applied to the above simple bilinear objective function, that is: 
\begin{equation}\label{eq: controlled simple GDA}
    z_{k+1}
    = (I -\eta M)z_k + u_k
\end{equation}
where $z_k=(x_k,y_k)\in\mathbb{R}^2$ represent the system state, $u_k\in\mathbb{R}^2$ a controlled input, and $M:=[\,0\,\,1\,;-1\,\,0\,]$.
In the absence of control ($u_k=0$, $\forall k$), the algorithm diverges as illustrated by the upper subplot of Figure \ref{fig:bilinear2d}.
To understand this unstable process, we keep track of the system's energy, expressed as the square 2-norm distance to the unique saddle point $(0,0)$, $V(z_k)=\|z_k\|^2=x_k^2+y_k^2$. A simple algebraic manipulation shows that 
\begin{align}\label{eq:V-controlled}
    V(z_{k+1}) &= (1+\eta^2)V(z_k) + 2u_k^Tz_k + \|u_k\|^2.
\end{align}
Notably, \eqref{eq:V-controlled} shows how, without control, the energy grows exponentially with each iteration.

We  to choose $u_k$ so that the right-hand side of \eqref{eq:V-controlled} is upper-bounded by $\alpha^2V(z_k)$, with $|\alpha|<1$. However, it is usually difficult to do without knowledge of the equilibrium.

The proposed approach uses state augmentation to incorporate a dissipation component that receives some external input $\hat u_k$ and filters its fluctuations, i.e., 
\begin{equation}\label{eq:dissipation}
\hat z_{k+1} = (1-\rho)\hat z_k + \rho \hat u_k. 
\end{equation}
The dissipation interpretation of \eqref{eq:dissipation} can be readily seen, again, by looking at the evolution of its energy
\begin{align}
&V(\hat{ z}_{k+1})
\!=\!(1\!-\!\rho)^2V(\hat z_k) \!+\!2\rho(1\!-\!\rho)\hat z_k^T\hat u_k\!+\!\rho^2\|\hat u_k\|^2\label{eq:V hat z}
\end{align}
which has the natural tendency to dissipate energy when $|1-\rho|<1$ and $\hat u_k = 0$, $\forall k$.

The key innovation of the proposed approach is to combine \eqref{eq: controlled simple GDA} and \eqref{eq:dissipation} with a proper choice of $u_k$ and $\hat u_k$ that satisfies:
\begin{itemize}
    \item The combined energy $V( \xi_k)=\|z_k\|^2+\|\hat z_k\|^2$ decreases, i.e., there is $|\alpha|<1$ such that:
    \begin{equation}\label{eq:V-linear decrease}
        V (\xi_{k+1})\leq\alpha^2 V(\xi_k)\,.
    \end{equation}
    \item The iterates $u_k\rightarrow u^*=0$, so that the equilibrium  of \eqref{eq: controlled simple GDA}, that is, the saddle point  does not change.
\end{itemize}
To achieve these goals we select
\begin{equation}\label{eq:u selection}
    u_k=-\rho (z_k-\hat z_k), \quad\text{ and }\quad
    \hat u_k =z_k \,,
\end{equation}
thus leading to the combined algorithm 
\begin{align}
    z_{k+1} & = z_k -\eta M z_k -\rho(z_k-\hat z_k)\\
    \hat z_{k+1} &= \hat z_k -\rho(\hat z_k - z_k)
\end{align}
The choice of $\hat u_k$ ensures that $\hat z_k$ acts as a filtering process of $z_k$. Further, since for constant  $z_k=\bar z$, $\hat z_k\rightarrow \bar z$, the choice of $u_k$ ensures that the equilibrium of \eqref{eq: controlled simple GDA} remains unchanged. We will show in Section \ref{Sec: Convergence} that one can choose $\eta$ and $\rho$ so that \eqref{eq:V-linear decrease} holds. We refer to the bottom of Figure \ref{fig:bilinear2d} for an illustration of its effectiveness.

\subsection{Dissipative GDA Algorithm}
In the previous section, we showed how to design a dissipative component \eqref{eq:dissipation} to stabilize the GDA updates on a simple bilinear objective. Since the design approach is independent of the particular saddle function $f(x,y)=xy$, it can be readily applied in more general settings. The algorithm presented below seeks to incorporate a similar friction component into a standard GDA update. This modification addresses the inherent limitation of the GDA update when applied to general bilinear problems.

\textbf{Dissipative gradient descent ascent (DGDA)}:
\begin{align}\label{eq: Dissipative}
    \begin{bmatrix}
     x_{k+1} \\
    \hat{ x}_{k+1}\\
     y_{k+1}\\
    \hat{ y}_{k+1}
    \end{bmatrix}= 
    \begin{bmatrix}
     x_k - \eta \nabla_ x f( x_k, y_k)- \rho( x_k-\hat{ x}_k) \\
    \hat{ x}_{k} - \rho(\hat{ x}_k -  x_k) \\
     y_k + \eta \nabla_ y f( x_k, y_k)- \rho( y_k-\hat{ y}_k) \\
    \hat{ y}_{k}-  \rho(\hat{ y}_k -  y_k)
    \end{bmatrix} 
\end{align}

Particularly, for $f$ as in \eqref{eq: min-max optimization problem}, in \eqref{eq: Dissipative} we introduce two new sets of variables $ \hat{ x} \in \mathbb{R}^n $ and $ \hat{ y} \in \mathbb{R}^m $ and a damping parameter $\rho > 0 $. 
One important observation is that, once the system reaches equilibrium, i.e., $x_{k+1} = x_k, y_{k+1} = y_k, \hat{ x}_{k+1} = \hat{ x}_{k} ,\hat{ y}_{k+1} = \hat{ y}_{k}  $, one necessarily has $\hat{ x}_k = x_k $ and $\hat{ y}_k = y_k $, which ensures that the fixed point is necessarily a critical point of the saddle function.

In the remaining part of this section, we will first provide several key observations and properties of the proposed DGDA updates. Then, in the next section, we will formally prove that the proposed DGDA algorithm provides a linear convergence guarantee for both bilinear and strongly convex-strongly concave functions. Mainly, we will show that by simply adding a friction or damping term, the DGDA updates succeed in dissipating the energy of the system even in cases where GDA does not.

\subsection{Key Properties and Related Algorithms}
The first important observation is that the above DGDA update could be considered as applying a vanilla GDA update to the following regularized surrogate for $f( x, y)$:
\begin{align}\label{eq:augmented saddle flow}
    f( x, y, \hat{ x}, \hat{ y})\!:=\! f( x, y) \!+\! \frac{\rho}{2} \| x \!-\!\hat{ x}  \|^2  \!-\! \frac{\rho}{2} \| y \!-\!\hat{ y}  \|^2.
\end{align}
Note that this is different from the \textit{Proximal Point Method} \cite{rockafellar1976monotone,parikh2014proximal} or introducing a \textit{$L_2$ regularization} \cite{krogh1991simple,hastie2009elements}.

\textbf{Differences with $L_2$ regularization}~
A commonly used method to ensure convergence is to introduce a $L_2$ regularization term in $x$ and $y$ \cite{krogh1991simple,hastie2009elements}:
\begin{align}
    \min_{ x \in \mathbb{R}^n}\max_{ y \in \mathbb{R}^m} f( x, y) +  \frac{\rho}{2} \| x\|^2 - \frac{\rho}{2} \| y \|^2.
\end{align}
Although the augmented objective function becomes strongly convex-strongly concave and the vanilla GDA updates will converge, this regularization changes the saddle points. While our algorithm also introduces two regularizing terms, the following Lemma verifies the fixed positions of saddle points between $ f( x, y)$ and $f( x, y, \hat{ x}, \hat{ y})$ with virtual variables aligned with original variables.    

\begin{lem}[Saddle Point Invariance]\cite[Lemma 6]{you2021saddle}
   For problem \ref{eq: min-max optimization problem}, a point $( x^*, y^*)$ is a saddle point of $ f( x, y)$ if and only if $( x^*, y^*,\hat{ x}^*,\hat{ y}^*)$ is a saddle point of  $f( x, y, \hat{ x}, \hat{ y}) $, with $\hat{ x}^* =  x^* $ and $\hat{ y}^* =  y^*$.
\end{lem}

More interestingly, the regularization terms, $\frac{\rho}{2} \| x -\hat{ x}  \|^2$ and $\frac{\rho}{2} \| y -\hat{ y}  \|^2$, do not introduce extra strong convexity-stong concavity to the original problem. Precisely, the augmented problem $ f( x, y, \hat{ x}, \hat{ y})$ is neither strongly convex on ($x,\hat{ x}$) nor strongly concave on $(y, \hat{ y})$. Indeed, on the hyperplane of $x = \hat{ x}$ and $y =\hat{ y} $, the augmented problem recovers the original problem $ f( x, y, \hat{ x}, \hat{ y}) = f( x, y)$. Additionally, the introduced regularization terms of the DGDA method are separable and local, thus preserving the distributed structure that original systems may have. Consequently, it can be seamlessly integrated into a fully distributed approach. 

\textbf{Differences with Proximal Point Method}
The Proximal Point Method for saddle point problems \cite{rockafellar1976monotone} shares a similar structure with DGDA algorithm. In the proximal method, the next iterates $(x_{k+1} ,y_{k+1})$ is the unique solution to the saddle point problem
\begin{equation}\label{eq:Proximal Point Lagrangian}
    \begin{aligned}
        & (x_{k+1},y_{k+1}) =   \mathrm{prox}_\eta(x_k,y_k)\\
         &\!:=\! \arg \min_{ x \in \mathbb{R}^n}\max_{ y \in \mathbb{R}^m} f( x, y) \!+\! \frac{\eta}{2}\|x\!-\!x_k\|^2 \!-\! \frac{\eta}{2}\|y\!-\!y_k\|^2
    \end{aligned}
\end{equation}
Using the optimality conditions of \eqref{eq:Proximal Point Lagrangian}, the update of the Proximal Point method can be written as:
\begin{equation}
    \begin{aligned}\label{eq: PPM}
   &x_{k+1} =  x_k - \eta \nabla_ x f( x_{k+1}, y_{k+1}), \\
     &y_{k+1} =  y_k + \eta \nabla_ y f( x_{k+1}, y_{k+1}).
\end{aligned}
\end{equation}
This expression shows that the Proximal Point method is an implicit method. Although Implicit methods are known to be more stable and to benefit from better convergence properties \cite{parikh2014proximal,gidel2018variational}, implementing the above updates requires computing the operators $(I+\eta \nabla_ x f)^{-1}$ and $(I+\eta \nabla_ y f)^{-1}$, and therefore may be computationally intractable. In contrast, the DGDA algorithm is an explicit algorithm, which applies a vanilla GDA update to the augmented objective function $\eqref{eq:augmented saddle flow}$. Thus, as shown in Table \ref{table: Comparison of bilinear convergence rate} and \ref{table: Comparison of str convergence rate}, and in the next section, we choose to compare the convergence  of DGDA only with other explicit algorithms for saddle point problems, such as 
the GDA, Extra-Gradient (EG), and Optimistic GDA (OGDA) methods which have comparable per-iteration computational requirements; although the  DGDA algorithm has twice as many state variables, it only requires a single gradient computation per update. Moreover, there is no need for retaining and reemploying the extrapolated gradient, which also sets it apart from the OGDA method. 

We finalize this section comparing DGDA with recent efforts leveraging \textit{Moreau-Yosida} smoothing techniques to solve nonconvex-concave \cite{zhang2020single,xu2023unified,yang2022faster}, nonconvex-nonconcave \cite{zheng2024universal} min-max optimization problems.

\textbf{Smoothed-GDA method} The \textit{Smoothed-GDA} consists of the following updates
\begin{equation}
    \begin{aligned}\label{eq: SGDA}
    & x_{k+1} =   x_k - c  \nabla_x K(x_k,z_k;y_k), \\
    & y_{k+1} =  y_k + \alpha \nabla_y K(x_{k+1},z_k;y_k)\\
    & z_{k+1} =   z_k + \beta(  x_{k+1} - z_k ) 
\end{aligned}
\end{equation}
where $  K( x, z;  y) = f(x,y) + \frac{p}{2} \|  x -  z\|^2.$
Smoothed-GDA was independently introduced by \cite{zhang2020proximal}. It was originally motivated as an ADMM algorithm to solve the linearly constrained nonconvex differentiable minimization problem \cite{zhang2020proximal}. The algorithm has been further extended to consider nonconvex-concave min-max optimization problems \cite{zhang2020single}. Despite its similarities with DGDA, the iterates of smoothed-GDA are sequential, requiring first updates in $x$ and subsequently in $(y,z)$, whereas DGDA updates are fully parallel or synchronized.
Unfortunately, the setting where such algorithm has been studied is different from the one considered in this paper, which difficult the comparison with the present work. A thorough comparison of DGDA with \eqref{eq: SGDA} is subject of current research.

\section{Convergence Analysis}\label{Sec: Convergence}
In this section, we provide a theoretical analysis of the proposed algorithm. 
\ifthenelse{\boolean{online}}{
}{
Due to space limitations, we have listed the key components of the proof but have omitted detailed calculations.
}
For the purpose of our analysis,  we consider a quadratic Lyapunov function to track the energy dissipation of the DGDA updates $$V_k:=   \| x_k - x^*\|^2 +   \| y_k - y^*\|^2 + \|\hat{ x}_k -\hat{ x}^* \|^2 + \| \hat{ y}_k -\hat{ y}^* \|^2,$$ which denotes the square 2-norm distance to the saddle point at the $k$-th iteration. 
The goal is, therefore, to find some $0 to \leq \alpha < 1$ such that $V_{k+1}\leq\alpha V_k$, where $\alpha$ denotes the linear convergence rate.

\subsection{Convergence Analysis for Bilinear Functions}
\label{subsec: Bilinear Convergence Analysis}
When applied to the bilinear min-max optimization problem $f(x,y) = x ^T A y$, the DGDA update \eqref{eq: Dissipative} is equivalent to a linear dynamical system. Specifically, denote $ z = [x,y]^T, \hat{z} = [\hat{x},\hat{y}]^T$ yields: 
\begin{align} \label{eq:Linear system bilinear functions}
    \begin{bmatrix}
    z_{k+1} \!-\! z^* \\
    \hat{z}_{k+1} \!-\! \hat{z}^*\\
    \end{bmatrix}
    \!=\! \begin{bmatrix}
    (1\!-\!\rho)I \!-\! \eta M \!&\!\!\! \rho I\\
      \rho I \!&\!\!\! (1\!-\!\rho)I \\
    \end{bmatrix} \!
    \begin{bmatrix}
    z_{k} \!-\! z^* \\
    \hat{z}_{k} \!-\! \hat{z}^*\\
    \end{bmatrix}, 
\end{align}
where $ M = \begin{bmatrix}
    \mathbf{0} & A \\
    -A^T & \mathbf{0}\\
    \end{bmatrix}. $
As a result, the linear convergence of DGDA, as well as its convergence rate can be derived from the analysis of the spectrum of the associated matrix that defines the DGDA update in \eqref{eq:Linear system bilinear functions}. This yields the following theorem.

\begin{thm} (Linear convergence of DGDA, Bilinear Case) \label{thm: Bilinear linear convergence}
    Let Assumption \ref{ass: Bilinear} hold. Then the updates \ref{eq: Dissipative} of DGDA with $0 < \eta \leq \frac{2\rho}{\sigma_{\max}(A)}$ and $\rho > 0$ provide linearly converging iterates. That is, there is a constant $\beta>0$, independent of $V_0$ such that
    \begin{align*}
     \textstyle  V_{k} \!\leq\! \mathcal{O}\left(\left( 1 \!-\! 2\rho \!+\! 2\rho^2 \!+\! (1\!-\!\rho)\sqrt{ 4\rho^2 \!-\!  \eta^2 \sigma_{\min}^2(A) }\right)^k \right) V_{0},
    \end{align*}
    Particularly, setting $\rho = 1/2$ and $\eta = 1/\sigma_{\max}(A)$  we have
    \begin{align}\label{eq: Convergence rate Bilinear}
\textstyle         V_{k} \leq \mathcal{O}\left(\left(1 -\frac{1}{4\kappa}\right)^k \right)\, V_{0}.
    \end{align}
    
\ifthenelse{\boolean{online}}{
\proof
The proof can be found in the appendix \ref{Appendix: Proof of Bilinear Theorem}
}{
\begin{proof} According to \cite[Lemma 7]{azizian2020tight}, we have $ \mathrm{Sp}(M) = \{ \pm i \sigma | \sigma^2 \in \mathrm{Sp}(AA^T) \}.$ Therefore, we can compute the eigenvalues of the system \eqref{eq:Linear system bilinear functions}: 
\begin{align}
    \mu_j =  1 - \rho  \pm i (\frac{1}{2} \eta \sigma_j) \pm \frac{1}{2} \sqrt{ 4\rho^2 - \eta^2 \sigma_j^2},
\end{align} where $\pm i\sigma_j \in  \mathrm{Sp}(M)$. 
Suppose that for all $j \in [m] $, we choose $0 < \eta \leq \frac{2\rho}{\sigma_{\max}} \leq \frac{2\rho}{\sigma_j}$ and $ \rho > 0$, which implies $4\rho^2 - \eta^2 \sigma_j^2 \geq 0$, then we can construct the following upper bound for the magnitude of eigenvalues,
\begin{align}
     |\mu_j|^2 & =  1 - 2\rho + 2\rho^2 \pm (1-\rho)\sqrt{ 4\rho^2 - \eta^2 \sigma_j^2 } \\
& <  1 - 2\rho + 2\rho^2 + (1-\rho) \sqrt{ 4\rho^2} \;\;=  1\;\; .
\end{align}
It follows from standard linear systems theory, e.g. \cite[Theorem 8.3]{hespanha2018linear}, the above spectral radius analysis of the linear system \eqref{eq:Linear system bilinear functions} results in the following linear convergence rate estimate:
\begin{align*}
   V_{k} \leq \mathcal{O}\left(\left( 1 - 2\rho + 2\rho^2 + (1-\rho)\sqrt{ 4\rho^2 -  \eta^2 \sigma_{\min}^2 } \right)^k \right) V_{0},
    \end{align*}
where $V_k:=   \| x_k - x^*\|^2 +   \| y_k - y^*\|^2 + \|\hat{ x}_k -\hat{ x}^* \|^2 + \| \hat{ y}_k -\hat{ y}^* \|^2$. Furthermore, the analysis of the above bound identifies the following optimal step sizes $ \eta = \frac{2\rho}{\sigma_{\max}}$ and $\rho = \frac{1}{2}$, and the following linear convergence rate estimate
\begin{align}
\textstyle         V_{k} \leq \mathcal{O}\left(\left(1 -\frac{1}{4\kappa}\right)^k \right)\, V_{0}.
\end{align}
\end{proof}
}
\end{thm}

We remark that linear convergence requires $\rho >0$. This is not surprising since GDA, which is known to diverge for bilinear functions, can be interpreted as DGDA method when $\rho = 0$. 
More importantly, by choosing the optimal step size $\rho=1/2 ,\eta=1/\sigma_{\max}(A 
)$, DGDA method achieves a better linear convergence rate than the EG and OGDA methods (see Table \ref{table: Comparison of bilinear convergence rate}). 

\subsection{Convergence Analysis for Strongly Convex Stronly Concave Functions}\label{subsec: Str Convergence Analysis}

We now consider the case of strongly convex-strongly concave min-max problems.
Let $F(z_k) := (\nabla_x f(x_{k},y_{k}), -\nabla_y f(x_{k},y_{k}))$. The DGDA updates can be written as follows:
\begin{align}
    \begin{bmatrix}
    z_{k+1} \\
    \hat{z}_{k+1}
    \end{bmatrix}= 
    \begin{bmatrix}
    z_k - \eta F(z_k)- \rho(z_k-\hat{z}_k) \\
    \hat{z}_{k} - \rho(\hat{z}_k -z_k)
    \end{bmatrix}
\end{align}

Because of the existence of the nonlinear term $F(z_k)$, we cannot analyze the spectrum as in the previous bilinear case. This is indeed a common challenge in analyzing most optimization algorithms beyond a neighborhood of the fixed point. 
We circumvent this problem by leveraging recent results on the analysis of variational mappings as $F(\cdot)$ via integral quadratic constraint~\cite{lessard2016analysis, hu2017dissipativity,fazlyab2018analysis}. 
\ifthenelse{\boolean{online}}{
More details can be found in the Appendix where we prove the following theorem.
}{
}

\begin{thm} (Linear convergence of DGDA, Strongly Convex-Strongly Concave Case) \label{thm: Strongly monotone linear convergence}
    Let Assumption \ref{ass: Str} hold, then the updates \eqref{eq: Dissipative} with $\rho = 1/2$ and $\eta = 1/(L + \mu) $ of the DGDA algorithm provide linearly converging iterates:
    \begin{align}\label{eq: Str Mon linear rate}
     V_{k} \leq \biggl (1 - \kappa^{-1} + \mathcal{O}( \kappa^{-2})
    \biggr) ^k  V_{0}
    \end{align} 
\ifthenelse{\boolean{online}}{
\proof
The proof can be found in Appendix \ref{Appendix: Proof of Strongly monotone Theorem}.
}{
\begin{proof}
Given a linear dynamical system of the form:
$ {\xi}_{k+1} = A {\xi} _{k} + B w_k, $
where ${\xi} \in \mathbb{R}^{n_{\xi}}$ is the state, $w_k \in \mathbb{R}^{n_{w}}$ is the input, $ A$ is the state transition matrix and $ B$ is the input matrix. Suppose that there exist a (Lyapunov) function $V$, satisfying $V({\xi}) \geq 0, \forall {\xi} \in \mathbb{R}^{n_{\xi}}$,  some $0 \leq \alpha < 1$ and a supply rate function $S( {\xi} _{k},w_k ) \leq 0 ,\forall k$ such that
\begin{align}\label{eq: dissipation inequality}
    V({\xi}_{k+1}) - \alpha^2  V({\xi}_{k}) \leq S( {\xi} _{k},w_k ), 
\end{align} then this dissipation inequality \eqref{eq: dissipation inequality} implies that $ V({\xi}_{k+1}) \leq \alpha^2  V({\xi}_{k}) $, and the state will approach a minimum value at equilibrium no slower than the linear rate $ \alpha^2$ \cite{hu2017dissipativity}. According to \cite[Lemma 6]{lessard2016analysis}, we could construct the following Linear Matrix Inequality and supply rate function for DGDA updates, by augmenting the states ${\xi} _{k} =(z_k ;\hat{z}_{k})$, \begin{align}
    S(\xi_k,w_k) \!=\! \begin{bmatrix}
         z_k \\
        \hat{z}_{k} \\
        w_k
    \end{bmatrix}^T\begin{bmatrix}
        2\mu LI & 0 &(\!-\!\mu\!+\!L)I \\
        0 & 0&0\\
        (\!-\!\mu\!+\!L)I & 0 & 2I \end{bmatrix}
        \begin{bmatrix}
         z_k \\
        \hat{z}_{k} \\
        w_k
    \end{bmatrix} \!\leq \! 0 
\end{align} where the nonlinear operator $F(z_k)$ meets the conditions specified in Assumption \ref{ass: Str}.

Finally, according to \cite[Theorem 2]{hu2017dissipativity}, constructing the dissipation inequality \eqref{eq: dissipation inequality} and proving linear convergence can be achieved through solving a semidefinite programming problem. Precisely, if there exists matrix $X^T = X$ and $P \in \mathbb{R}^{n_{\xi} \times n_{\xi} }$ with $P \succeq 0$ such that 
\begin{align}\label{eq: Dis LMI}
    \begin{bmatrix}
        A^T P A - \alpha^2 P & A^T P B \\
        B^T P A & B^T P B
    \end{bmatrix} -X \leq 0,
\end{align}
where $ S( {\xi} ,w ) :=\begin{bmatrix}
            {\xi} \\ w
        \end{bmatrix}^T X\begin{bmatrix}
            {\xi} \\ w
        \end{bmatrix}, $
then the dissipation inequality holds for all trajectories of $ {\xi}_{k+1} = A {\xi} _{k} + B w_k, $ with $ V({\xi}) ={\xi}^T P {\xi} $. Given the set of problem parameters,
a set of feasible solutions is given by: \begin{align}
 \rho = \frac{1}{2} ,\eta = \frac{1}{L+\mu}, 
  P = \begin{bmatrix}
      (L+\mu)^2 & 0\\
      0 & (L+\mu)^2
  \end{bmatrix}\otimes I, \\
  \alpha^2 = \frac{3L^2 + 2L \mu + 3\mu^2 + \sqrt{(L+\mu)^4 + 16L^2\mu^2}}{4(L+\mu)^2} .
\end{align}

After substituting the condition number $\kappa:=L /\mu$, the convergence rate simplifies to $  \alpha^2 = 1 - \kappa^{-1} + \mathcal{O}\bigl((\frac{\mu}{L})^2 \bigr) $
\end{proof}
}
\end{thm}

Similarly, as in the bilinear case, we remark on the importance of the dissipation component. When $\rho =0$, a similar analysis as in the proof of the theorem recovers the lower bound of the convergence rate of GDA $(1 - \kappa^{-2})$ as shown in \cite[3.1]{zhang2021unified}. Thus, our DGDA method provides a better convergence rate estimate than GDA, since clearly $ \kappa \in [1,\infty)$, and therefore $\kappa^{-2} \leq \kappa^{-1}$.

We remark that while the rate obtained in Theorem \ref{thm: Strongly monotone linear convergence} is clearly better than those of the EG and OGDA methods for large condition numbers $\kappa$ (see Table \ref{table: Comparison of str convergence rate}), the theorem fails to quantify the comparative performance of DGDA for small values of $\kappa$.
The following corollary shows that indeed, the rate of DGDA is provably better for all $\kappa\geq 2$.


\begin{cor}[SCSC, comparison with known rates] \label{cor:Str montonte comparison with known rates}
Let Assumption \ref{ass: Str} hold, and suppose that $L \geq 2m$, i.e., $\kappa \geq 2$. Then, the linear convergence rate estimate of DGDA \eqref{eq: Str Mon linear rate} is smaller (better) than that of of EG and OGDA, i.e.,  $1 - \kappa^{-1}/4 $ (Theorem $6 \& 7$ \cite{azizian2020tight} and Theorem $4\&7$ \cite{mokhtari2020unified}).

\end{cor}

\section{Numerical Experiments}\label{Sec: Numerical Experiments}
\begin{figure}[t]
\centering
\includegraphics[width=0.8\columnwidth]{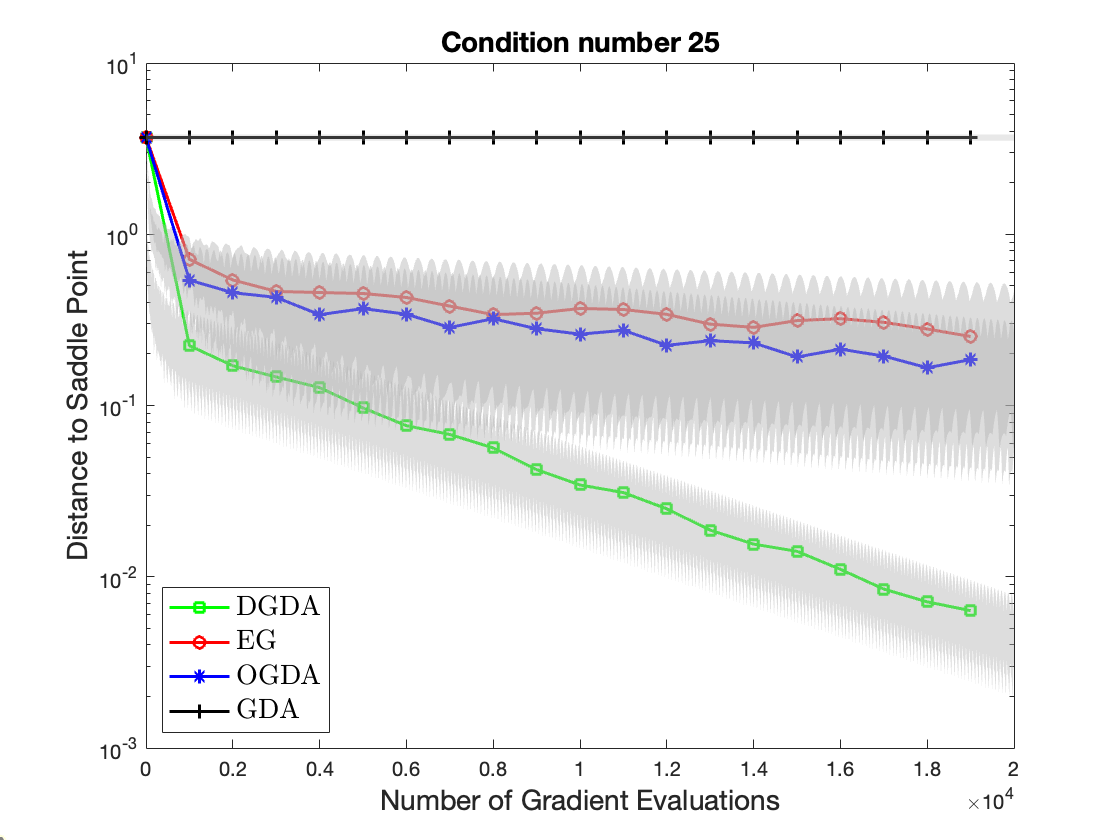}
    \caption{Convergence of GDA, EG, OGDA, and DGDA in terms of the number of gradient evaluations for the bilinear problem. GDA diverges and the error is not shown.
    All other three algorithms converge linearly, where the DGDA method provides the best performance.  }
    \label{fig:distance_bilinear}
\end{figure}

In this section, we compare the performance of the proposed Dissipative gradient descent (DGDA) method with the Extra-gradient (EG), Gradient descent ascent (GDA), and Optimistic gradient descent ascent (OGDA) methods.

\subsection{Bilinear problem}\label{sec: numerical bilinear}

\begin{figure}[h]
\centering
\includegraphics[width=0.8\columnwidth]{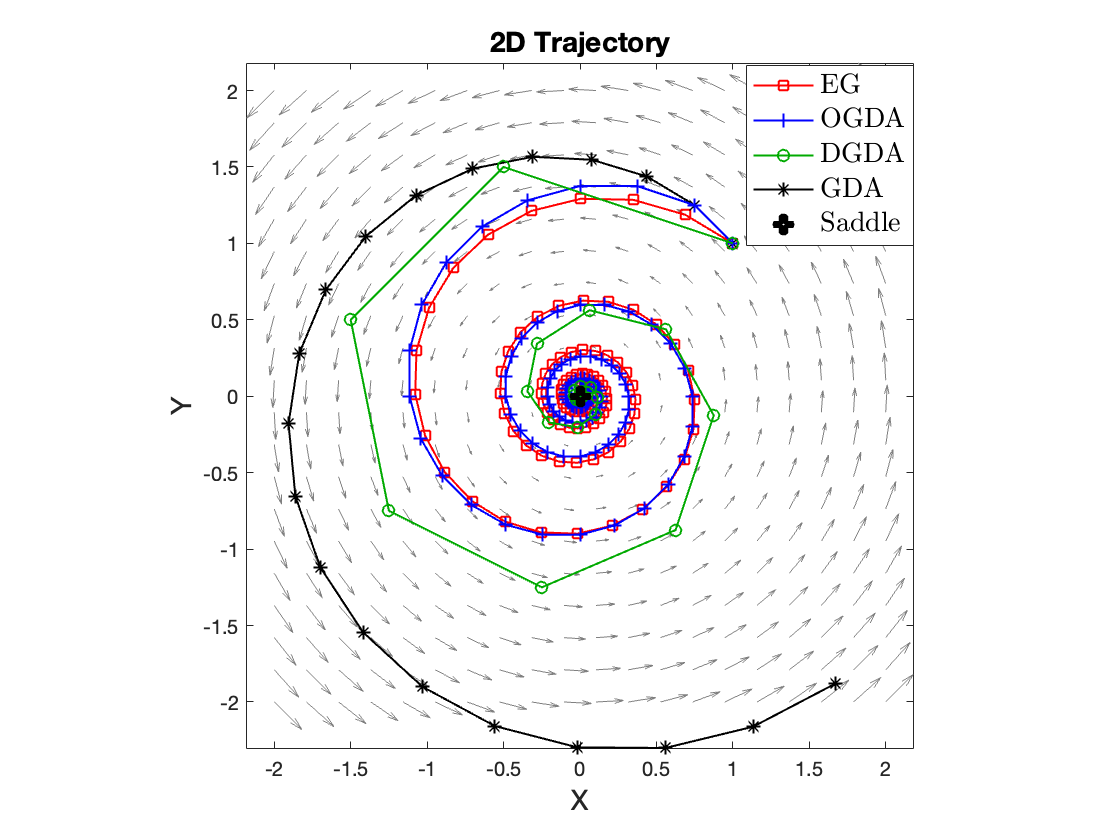}
    \caption{Trajectories of GDA, EG, OGDA, and DGDA for a 2d bilinear problem.  GDA diverges and all other three algorithms converge linearly, where the DGDA method provides the best performance.    }
    \label{fig:2d}
\end{figure}

We first consider the following bilinear min-max optimization problem: $\min_{x \in \mathbb{R}^n}\max_{y \in \mathbb{R}^m} x^T A y$,  where $A \in \mathbb{R}^{m \times n}$ is full-rank. The simulation results are illustrated in Figure \ref{fig:distance_bilinear} and Figure \ref{fig:2d}. In this experiment, we set the dimension of the problem to $m = n = 10$ and the iterates are initialized at $x_0, y_0 $, which are randomly drawn from the uniform distribution on the open interval $(0,1)$. 

We plot the errors (distance to saddle points) of DGDA, EG, and OGDA versus the number of gradient evaluations for this problem in Figure \ref{fig:distance_bilinear}. The solid line and grey-shaded error bars represent the average trajectories and standard deviations of 20 trials, where in each trial the randomly generated matrix  $A$ has a fixed condition number, i.e., $\kappa =  \sigma^2_{\max}(A)/ \sigma^2_{\min}(A) = 25$. The key motivation is that all three algorithms' convergence rates critically depend on $\kappa^{-1}$, and by fixing the condition number, we provide an explicit comparison of their convergence speed. 

We pick the step size for different methods according to theoretical findings. That is, we select  $\rho = 1/2$ and $\eta = 1/\sigma_{\max}(A)$ for DGDA (Theorem \ref{thm: Bilinear linear convergence}), $\eta = 1/4L = 1/4\sigma_{\max}(A)$ for EG and OGDA (Theorem 6\&7 \cite{azizian2020tight} and Theorem 4\&7 \cite{mokhtari2020unified}). In Figure \ref{fig:distance_bilinear}, we do not show the error of GDA since it diverges for this 
bilinear saddle point problem. All other three algorithms converge linearly, with the DGDA method providing the best performance. 

Finally, to provide a qualitative demonstration of how DGDA fares with other existing algorithms, we further plot the sample trajectories of GDA, EG, OGDA, and EGDA on a simple 2D bilinear min-max problem, with $m = n = 1$. In Figure \ref{fig:2d}, we observe that while GDA diverges, the trajectories of all other three algorithms converge linearly to the saddle point $(x^*, y^*) = (0,0)$. Interestingly, our proposed algorithm (DGDA) despite taking larger steps, exhibits faster linear convergence.
\subsection{Strongly convex-strongly concave problem} 
\begin{figure}[h]
\centering
\includegraphics[width=0.9\columnwidth]{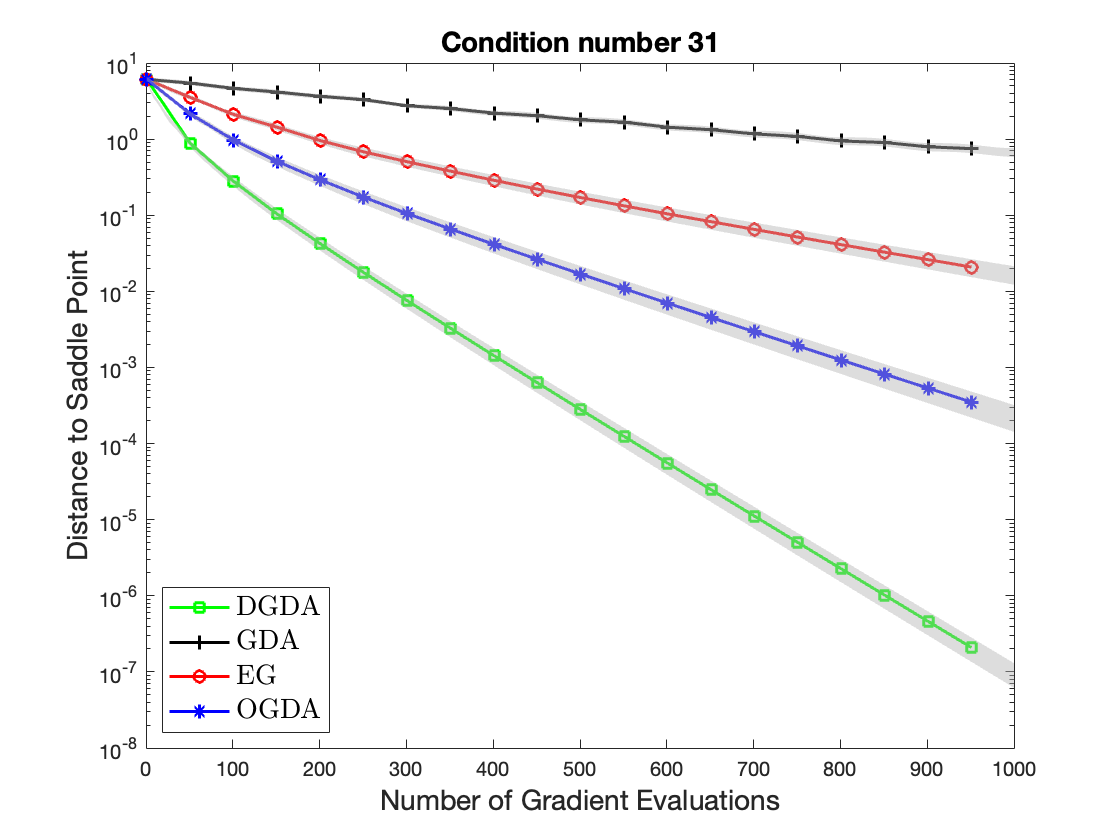}
    \caption{Convergence of GDA, EG, OGDA, and DGDA in terms of the number of gradient evaluations for problem \ref{eq: Numerical Str}. All algorithms converge linearly, and the DGDA method has the best performance.  }
    \label{fig:str}
\end{figure}

In the second numerical example, we focus on a strongly convex-strongly strongly-concave quadratic problem of the following form:
\begin{align}\label{eq: Numerical Str}
      \min_{x \in \mathbb{R}^n}\max_{y \in \mathbb{R}^m} \frac{1}{2} x^T A x - \frac{1}{2} y^T B y +  x^T C y,
\end{align}
where the matrices satisfy $\mu_A I \preceq A \preceq L_A I $,  $\mu_B I \preceq B \preceq L_B I $, $\mu_c^2 I \preceq C^T C \preceq L_c^2 I $. As a result, the problem \eqref{eq: Numerical Str} satisfy Assumption \ref{ass: Str}. In this experiment, we set the dimension of the problem to $n = 50, m = 10$, and the iterates are initialized at $x_0, y_0$, which are randomly drawn from the uniform distribution on the open interval $(0,1)$. We plot the errors (distance to saddle points) of GDA, DGDA, EG, and OGDA versus the number of gradient evaluations for this problem in Figure \ref{fig:str}. Again, the solid line and grey-shaded error bars represent the average trajectories and standard deviations of 20 trials, where in each trial the randomly generated matrix $
 \begin{bmatrix}
    A & C \\
    - C^T & B
\end{bmatrix}    
$
is chosen such that the condition number of \eqref{eq: Numerical Str} remains constant, i.e., $\kappa =  L/\mu= 31$. Similarly as in the bilinear problem in Section \ref{sec: numerical bilinear}, we pick the step size for the DGDA method according to our theoretical finding in Theorem \ref{thm: Strongly monotone linear convergence}.  The step size of the GDA method is selected as $\eta = \mu/L^2$ (Theorem 5 \cite{beznosikov2023smooth}). The step sizes for EG and OGDA methods are selected as $\eta = 1/4L $  (Theorem 6\&7 \cite{azizian2020tight} and Theorem 4\&7 \cite{mokhtari2020unified}). According to the plots, all algorithms converge linearly, and the DGDA method has the best performance.

\section{Conclusion and Future Work}\label{Sec: Conclusion}
In this work, we present the Dissipative GDA (DGDA) algorithm, a novel method for solving min-max optimization problems. Drawing inspiration from dissipativity theory and control theory, we address the challenge of diverging oscillations in bilinear min-max optimization problems when using the Gradient Descent Ascent (GDA) method. Particularly, we introduce a friction term into the GDA updates aiming to dissipate the internal energy and drive the system towards equilibrium. By incorporating a state-augmented regularization, our proposed DGDA method can be seen as performing standard GDA on an extended saddle function without introducing additional convexity. We further establish the superiority of the convergence rate of the proposed DGDA method when compared with other established methods including GDA, Extra-Gradient (EG), and Optimistic GDA. The analysis is further supported by two numerical examples, demonstrating its effectiveness in solving saddle point problems. Our future work includes studying the DGDA method in a stochastic setting and its application in solving Constrained Reinforcement learning problems in the policy space.



\bibliographystyle{IEEEtran}
\bibliography{main.bib} 



\appendix
\subsection{First Order Algorithms for Saddle Point Problems} \label{Appendix: First Order Algorithms}

In this section, we introduce several popular first-order methods for solving the min-max problem \ref{eq: min-max optimization problem} in the machine learning community. Precisely, we focus on Gradient Descent-Ascent (GDA), Extra-gradient (EG), Optimistic Gradient Descent-Ascent (OGDA), Proximal Point and Smoothed-GDA methods.

\subsubsection*{Gradient descent ascent (GDA)}
\begin{equation}
\begin{aligned}
     &x_{k+1} =  x_k - \eta \nabla_ x f( x_k, y_k), \\
     &y_{k+1} =  y_k + \eta \nabla_ y f( x_k, y_k).
\end{aligned}    
\end{equation}

When the problem is strongly convex-strongly concave and L-Lipschitz, the GDA method provides linear convergence, with step size $\eta = \mu/L^2$ and a know rate estimate of $ 1 - 1/\kappa^2$ \cite{grimmer2023landscape,zhang2021unified}. However, when the problem is bilinear, the standard GDA method fails to converge. Therefore, variants of the gradient method such as Extra-Gradient and Optimistic Gradient Descent-Ascent methods have attracted much attention in recent literature because of their superior empirical performance in solving min-max optimization problems such as training GANs and solving C-RL problems.

\subsubsection*{Extra-gradient (EG)}

\begin{equation}
    \begin{aligned}\label{eq: EG}
    & x_{k+1/2} =  x_k - \eta \nabla_ x f( x_{k}, y_{k}), \\
    & y_{k+1/2} =  y_k + \eta \nabla_ y f( x_{k}, y_{k}),  \\
    & x_{k+1} =  x_k - \eta \nabla_ x f( x_{k+1/2}, y_{k+1/2}), \\
    & y_{k+1} =  y_k + \eta \nabla_ y f( x_{k+1/2}, y_{k+1/2}).
\end{aligned}
\end{equation}

Extra-gradient is a classical method introduced in \cite{korpelevich1976extragradient}, where its linear rate of convergence for bilinear functions and smooth strongly convex-strongly concave functions have been established in many recent works (see Table \ref{table: Comparison of bilinear convergence rate} and \ref{table: Comparison of str convergence rate}). The Extra-gradient method first computes an extrapolated point ($ x_{k+1/2},  y_{k+1/2}$) by performing a GDA update. Then the gradients evaluated at the extrapolated point are used to compute the new iterates $( x_{k+1}, y_{k+1}) $.

The linear convergence rate of EG for strongly convex-strongly concave is established, with a standard known rate of $ 1 - 1/4\kappa$; see e.g.\cite{azizian2020tight,mokhtari2020unified}. One issue with the Extra-gradient method is that, as the name suggests, each update requires evaluation of extra gradients at the extrapolated point ($ x_{k+1/2},  y_{k+1/2}$), which doubles the computational complexity of EG method compared to vanilla GDA method. 

\subsubsection*{Optimistic gradient descent ascent (OGDA)}

\begin{equation}
    \begin{aligned}\label{eq: OGDA}
    & x_{k+1} \! =\!  x_k - 2\eta \nabla_ x f( x_k, y_k) +  \eta \nabla_ x f( x_{k-1}, y_{k-1}), \\
    & y_{k+1} \!= \! y_k + 2\eta \nabla_ y f( x_k, y_k) - \eta \nabla_ y f( x_{k-1}, y_{k-1}).
\end{aligned}
\end{equation}

The Optimistic gradient descent ascent (OGDA) method adds a "negative-momentum" term to each of the updates, which differentiates the OGDA method from the vanilla GDA method. Meanwhile, the OGDA method stores and re-uses the extrapolated gradient for the extrapolation, which only requires a single gradient computation per update.

The convergence properties of OGDA were also recently investigated in (refer to Table \ref{table: Comparison of bilinear convergence rate} and \ref{table: Comparison of str convergence rate}), demonstrating linear convergence rates with smooth and bilinear functions, as well as strongly convex-strongly concave functions.

\subsubsection*{Proximal Point (PP)}
The proximal point method for convex minimization has been extensively studied \cite{rockafellar1976augmented,parikh2014proximal} and extended to solve saddle point problems in \cite{rockafellar1976monotone}. Define the iterates $\{ x_{k+1} ,y_{k+1}\}$ as the unique solution to the saddle point problem
\begin{equation}
    \begin{aligned}
        & (x_{k+1},y_{k+1}) =   \mathrm{prox}_\eta(x_k,y_k)\\
         &:= \arg \min_{ x \in \mathbb{R}^n}\max_{ y \in \mathbb{R}^m} f( x, y) + \frac{\eta}{2}\|x-x_k\|^2 - \frac{\eta}{2}\|y-y_k\|^2
    \end{aligned}
\end{equation}
Using the optimality conditions, the update of the Proximal Point method can be written as:
\begin{equation}
    \begin{aligned}
   &x_{k+1} =  x_k - \eta \nabla_ x f( x_{k+1}, y_{k+1}), \\
     &y_{k+1} =  y_k + \eta \nabla_ y f( x_{k+1}, y_{k+1}).
\end{aligned}
\end{equation}
This expression shows that in contrast to explicit methods such as GDA, EG, and OGDA methods, the Proximal Point method is an implicit method. Although Implicit methods are known to be more stable and to benefit from better convergence properties \cite{parikh2014proximal,gidel2018variational}, implementing the above updates requires computing the operators $(I+\eta \nabla_ x f)^{-1}$ and $(I+\eta \nabla_ y f)^{-1}$, and therefore may be computationally intractable. Notably, in \cite{mokhtari2020unified}, the authors show the EG and OGDA methods can be interpreted as an approximation of the PP method, and therefore exhibits similar convergence performance.

\subsubsection*{Smoothed-GDA method}
We finalize this section comparing DGDA with recent efforts leveraging \textit{Moreau-Yosida} smoothing techniques to solve nonconvex-concave \cite{zhang2020single,xu2023unified,yang2022faster}, nonconvex-nonconcave \cite{zheng2024universal} min-max optimization problems.
\begin{equation}
    \begin{aligned}
    & x_{k+1} =   x_k - c  \nabla_x K(x_k,z_k;y_k), \\
    & y_{k+1} =  y_k + \alpha \nabla_ K(x_{k+1},z_k;y_k)\\
    & z_{k+1} =   z_k + \beta(  x_{k+1} - z_k ) 
\end{aligned}
\end{equation}
where
    $K( x, z;  y) = f(x,y) + \frac{p}{2} \|  x -  z\|^2$.

The Smoothed-GDA was independently introduced by Jiawei et al. in \cite{zhang2020proximal} and later \cite{zhang2020single}. It was originally motivated by ADMM to solve the linearly constrained nonconvex differentiable minimization problem \cite{zhang2020proximal}, where they introduce an extra quadratic proximal term for the equality constraints and an extra sequence $\{ z_k\}$. They claim this smoothing or exponential averaging scheme is necessary for the convergence of the proximal ADMM when the objective is nonconvex. Later on, this scheme is further extended to solve the nonconvex-concave min-max optimization problem \cite{zhang2020single}.

\subsection{Proof of Theorem \ref{thm: Bilinear linear convergence}} \label{Appendix: Proof of Bilinear Theorem}
We consider, for ease of presentation, the case when $A \in \mathbb{R}^{m \times m}$ is a square matrix. The extension for non-square matrices is straightforward and has been covered in the literature \cite[Appendix G]{zhang2019convergence}. Applying the updates \ref{eq: Dissipative} to $f(x,y) = x ^T A y$ and denoting $ z = [x,y]^T, \hat{z} = [\hat{x},\hat{y}]^T$ yields: 
\begin{equation}\label{eq: linear system bilinear DGDA updates}
    \begin{aligned}
         \begin{bmatrix}
    z_{k+1} - z^* \\
    \hat{z}_{k+1} - \hat{z}^*\\
    \end{bmatrix} &= 
    \begin{bmatrix}
    z_k - \eta M z_k- \rho(z_k-\hat{z}_k) \\
      \hat{z}_{k}-  \rho(\hat{z}_k- z_{k}) \\
    \end{bmatrix} \\
    &= \begin{bmatrix}
    (1-\rho)I - \eta M & \rho I\\
      \rho I & (1-\rho)I \\
    \end{bmatrix} 
    \begin{bmatrix}
    z_{k} - z^* \\
    \hat{z}_{k} - \hat{z}^*\\
    \end{bmatrix},
    \end{aligned}
\end{equation}
where 
\begin{align*}
    M = \begin{bmatrix}
    \mathbf{0} & A \\
    -A^T & \mathbf{0}\\
    \end{bmatrix} 
\end{align*}
According to \cite[Lemma 7]{azizian2020tight}
\begin{align*}
    \mathrm{Sp}(M) = \{ \pm i \sigma | \sigma^2 \in \mathrm{Sp}(AA^T) \}.
\end{align*}
We will use $\sigma_{\max} $ and $\sigma_{\min}$ to denote the largest singular value and smallest singular of matrix $A$, respectively. And according to Assumption \ref{ass: Bilinear}, we have $\sigma_{\min} > 0$.
Since $M$ is a normal matrix and diagonalizable, we can compute the eigenvalues of the linear system \eqref{eq: linear system bilinear DGDA updates} using the following similarity transformation
\begin{equation}\label{eq: Spectral Analysis Matrix for Bilinear}
    \begin{aligned}
       & \begin{bmatrix}
    (1-\rho)I - \eta M & \rho I\\
      \rho I & (1-\rho)I \\
    \end{bmatrix}  = \\ &\begin{bmatrix}
    U^{-1}& \mathbf{0}\\
      \mathbf{0} & U^{-1} \\
    \end{bmatrix}
    \begin{bmatrix}
    (1-\rho) I - \eta \boldsymbol\Lambda &\rho I \\
      \rho I  & (1-\rho) I \\
    \end{bmatrix}
    \begin{bmatrix}
    U& \mathbf{0}\\
      \mathbf{0} & U \\
    \end{bmatrix},
    \end{aligned}
\end{equation}

where $\boldsymbol\Lambda = \mathrm{diag}(\lambda_1,...,\lambda_{2m}) $, $\lambda_{2j-1} = + i \sigma_j$, $\lambda_{2j} = - i \sigma_j$, with $\pm i\sigma_j \in  \mathrm{Sp}(M), j = \{1, ...,m \}$. In order to show linear convergence, we want to show that $ \max_{j \in [m]}{|\mu_j|^2} < 1$ , where $\mu_j$ are the eigenvalues of the above matrix \eqref{eq: Spectral Analysis Matrix for Bilinear} (i.e., the spectral radius of matrix  \eqref{eq: Spectral Analysis Matrix for Bilinear}). As straight forward computation leads to
\begin{align}
    \mu_j = &\frac{1}{2} (2- 2\rho -\eta \lambda_j \pm \sqrt{ \eta^2 \lambda_j^2 + 4\rho^2} ) \nonumber \\
    = & 1 - \rho  \pm i (\frac{1}{2} \eta \sigma_j) \pm \frac{1}{2} \sqrt{ 4\rho^2 - \eta^2 \sigma_j^2}
\end{align}
Since, for complex number $c$, $|c|^2  = c \Bar{c} $,  the magnitude of eigenvalues  $ |\mu_j|^2$ are given by,
\begin{equation}\label{eq:bilinear eigenvalue magnitue}
    \begin{aligned}
        &|\mu_j|^2 = \biggl ( 1 - \rho  + i \frac{1}{2} \eta \sigma_j \pm \frac{1}{2} \sqrt{ 4\rho^2 - \eta^2 \sigma_j^2} \biggr) \times \\ 
        &\qquad \quad \; \biggl ( 1 - \rho  - i \frac{1}{2} \eta \sigma_j \pm \frac{1}{2} \overline{\sqrt{ 4\rho^2 - \eta^2 \sigma_j^2}} \biggr ) \nonumber\\
    &= (1\!-\!\rho)^2 \!+\! \frac{1}{4}\eta^2 \sigma_j^2 \!+\! \frac{1}{4}|4\rho^2 \!-\! \eta^2 \sigma_j^2   \!\pm \!(1\!-\!\rho)\Re(\sqrt{ 4\rho^2 \!-\! \eta^2 \sigma_j^2}) \nonumber\\
     &\quad \pm \frac{i}{4}\eta \sigma_j\overline{\sqrt{ 4\rho^2 - \eta^2 \sigma_j^2}} \mp \frac{i}{4}\eta \sigma_j\sqrt{ 4\rho^2 - \eta^2 \sigma_j^2} \nonumber \\
     &=  \begin{cases}
     1 \!- \!2\rho \!+ \!2\rho^2  \!\pm \!(1\!-\!\rho)\sqrt{ 4\rho^2 \!- \!\eta^2 \sigma_j^2} , \;\text{if $4\rho^2 \!-\! \eta^2 \sigma_j^2 \!\geq \!0 $}\\
      1 \!-\! 2\rho \!+ \!\frac{1}{2}\eta^2 \sigma_j^2  \! \pm \!\frac{1}{2}\eta \sigma_j \sqrt{ \eta^2 \sigma_j^2 \!-\! 4\rho^2 }, \;\text{if $ \eta^2 \sigma_j^2 \!-\! 4\rho^2 \!\geq\! 0 $}
		 \end{cases}
    \end{aligned}
\end{equation}

Suppose that for all $j \in [m] $, we choose $0 < \eta \leq \frac{2\rho}{\sigma_{\max}} \leq \frac{2\rho}{\sigma_j}$ and $ \rho > 0$, which implies $4\rho^2 - \eta^2 \sigma_j^2 \geq 0$. From \eqref{eq:bilinear eigenvalue magnitue}, $\forall j \in [m]$ we have,
\begin{align*}
|\mu_j|^2 & =  1 - 2\rho + 2\rho^2 \pm (1-\rho)\sqrt{ 4\rho^2 - \eta^2 \sigma_j^2 } \\
& \leq  1 - 2\rho + 2\rho^2 + (1-\rho)\sqrt{ 4\rho^2 -  \eta^2 \sigma_{j}^2 } \\
& \leq  1 - 2\rho + 2\rho^2 + (1-\rho)\sqrt{ 4\rho^2 -  \eta^2 \sigma_{\min}^2 } \\
& <  1 - 2\rho + 2\rho^2 + (1-\rho) \sqrt{ 4\rho^2} \\
& = 1\;\; .
\end{align*}
According to classical linear system theory, e.g. \cite[Theorem 8.3]{hespanha2018linear}, the above spectral radius analysis of the linear system \eqref{eq: linear system bilinear DGDA updates} results in the following linear convergence rate estimate:
\begin{align*}
   V_{k} \leq \mathcal{O}\left(\!\left( 1 - 2\rho + 2\rho^2 + (1-\rho)\sqrt{ 4\rho^2 -  \eta^2 \sigma_{\min}^2 } \right)^k \!\right) V_{0},
    \end{align*}
where $V_k:=   \| x_k - x^*\|^2 +   \| y_k - y^*\|^2 + \|\hat{ x}_k -\hat{ x}^* \|^2 + \| \hat{ y}_k -\hat{ y}^* \|^2$.

Furthermore, we want to select the optimal step size $\rho, \eta$. The immediate step is to substitute the optimal $ \eta = \frac{2\rho}{\sigma_{\max}}$, which yields the following inequality:
\begin{align}
    |\mu_j|^2 & \leq  1 \!- \!2\rho \!+ \!2\rho^2 \!+ \!(1-\rho)\sqrt{ 4\rho^2 - \frac{4\rho^2}{\sigma_{\max}^2} \sigma_j^2 } \;\;, \forall j \in [m].
\end{align}
The spectral radius is therefore given by choosing $\sigma_j = \sigma_{\min} $ above, i.e.,
\begin{align*}
  \max_{j \in [m]}&{|\mu_j|^2} =  2\rho^2 - 2\rho + 1 +(1-\rho)2\rho \sqrt{1 - \frac{\sigma_{\min}^2}{\sigma_{\max}^2}  } \\
 & = \rho^2 \biggl (2\!-\!2\sqrt{1 - \frac{\sigma_{\min}^2}{\sigma_{\max}^2}  } \biggr) \!- \!\rho \biggl (2\!-\!2\sqrt{1 - \frac{\sigma_{\min}^2}{\sigma_{\max}^2}  }\biggr ) \!+ \!1 \\
    & \leq 1 - \frac{1}{2}(1-\sqrt{1 - \frac{\sigma_{\min}^2}{\sigma_{\max}^2}  }) \\
    &= \frac{1}{2} + \frac{1}{2}\sqrt{1 - \frac{\sigma_{\min}^2}{\sigma_{\max}^2}  } \;\;, \forall j \in [m]
\end{align*}
where the last inequality comes from selecting optimal $\rho = \frac{1}{2}$ of a quadratic polynomial  of $\rho$.
Using the fact that $\sqrt{1-x} \leq 1- x/2 $, we have 
\begin{align}
    \max_{j \in [m]}{|\mu_j|^2} \leq 1 -\frac{1}{4} \frac{\sigma_{\min}^2}{\sigma_{\max}^2} 
\end{align}
Again, this results in the following linear convergence rate estimate:
\begin{align}
\textstyle         V_{k} \leq \mathcal{O}\left(\left(1 -\frac{1}{4\kappa}\right)^k \right)\, V_{0}.
\end{align}

\begin{rem}
We could also choose $\eta = \frac{2\rho}{\sigma_{\min}} $ such that $ \eta^2 \sigma_j^2 - 4\rho^2 \geq 0 $. And we could construct a similar linear convergence rate by repeating the above process. However, in practice, we found that the step sizes $\eta = \frac{2\rho}{\sigma_{\max}}, \rho = 1/2 $  always perform better in numerical experiments. Therefore, we choose this pair of step sizes by default. 
\end{rem} 

\begin{rem}
Since GDA method could be interpreted as a special case of DGDA method when selecting $\rho = 0$, the above step proves that when $ \eta > 0$, the GDA method diverges for a bilinear objective function. Specifically, when $\rho = 0 ,\eta > 0 $, we have $ \eta^2 \sigma_j^2- 4\rho^2 > 0 $ and 
\begin{align}
    |\mu_j|^2 = 1 \pm \frac{1}{2}\eta^2 \sigma^2_j,
\end{align}
\end{rem}


\subsection{Proof of Theorem \ref{thm: Strongly monotone linear convergence}} \label{Appendix: Proof of Strongly monotone Theorem}
The proof relies on the application of dissipativity theory to construct Lyapunov functions and establish linear convergence. For more detailed information, refer to \cite{hu2017dissipativity}.

According to \cite{hu2017dissipativity}, a linear dynamical system of the form:
\begin{align}\label{eq: dis dynamical system}
  {\xi}_{k+1} = A {\xi} _{k} + B w_k
\end{align}
Here, ${\xi} \in \mathbb{R}^{n_{\xi}}$ is the state, $w_k \in \mathbb{R}^{n_{w}}$ is the input, $ A$ is the state transition matrix and $ B$ is the input matrix. Suppose that there exist a (Lyapunov) function $V$, satisfying $V({\xi}) \geq 0, \forall {\xi} \in \mathbb{R}^{n_{\xi}}$,  some $0 \leq \alpha < 1$ and a supply rate function $S( {\xi} _{k},w_k ) \leq 0 ,\forall k$ such that
\begin{align}\label{eq: dissipation inequality}
    V({\xi}_{k+1}) - \alpha^2  V({\xi}_{k}) \leq S( {\xi} _{k},w_k ).
\end{align}
This dissipation inequality \eqref{eq: dissipation inequality} implies that $ V({\xi}_{k+1}) \leq \alpha^2  V({\xi}_{k}) $, and the state will approach a minimum value ate equilibrium no slower than the linear rate $ \alpha^2$.
The flowing theorem states how to construct the dissipation inequality \eqref{eq: dissipation inequality} by solving a semidefinite programming problem.
\begin{thm}\cite{hu2017dissipativity}[Theorem 2] \label{thm: Dissipative}
    Consider the following quadratic supply rate with $X \in \mathbb{R}^{(n_{\xi} + n_{w}) \times (n_{\xi} + n_{w})}$ and $X^T = X $
    \begin{align}
        S( {\xi} ,w ) :=\begin{bmatrix}
            {\xi} \\ w
        \end{bmatrix}^T X\begin{bmatrix}
            {\xi} \\ w
        \end{bmatrix}.
    \end{align}
If there exists matrix $P \in \mathbb{R}^{n_{\xi} \times n_{\xi} }$ with $P \succeq 0$ such that 
\begin{align}\label{eq: Dis LMI}
    \begin{bmatrix}
        A^T P A - \alpha^2 P & A^T P B \\
        B^T P A & B^T P B
    \end{bmatrix} -X \leq 0,
\end{align}
     then the dissipation inequality  holds for all trajectories of \eqref{eq: dis dynamical system} with $ V({\xi}) ={\xi}^T P {\xi} $.
\end{thm}
A major benefit of the proposed constructive dissipation approach is that it replaces the trouble some component of a dynamical system (e.g. the gradient term $w=\nabla_{{\xi}} f({\xi})$) by a quadratic constraint on its inputs and outputs that is always satisfied, namely the supply rate constraint $S(\xi,w)\leq 0$. This leads to a two-step novel approach to the convergence analysis of optimization algorithms.
\begin{enumerate}
    \item Choose a proper quadratic supply rate function $S$ such that $S( {\xi} _{k},w_k ) \leq 0 ,\forall k$, that depends on the specific nonlinear term.
    \item Solve the Linear Matrix Inequality \eqref{eq: Dis LMI} to obtain a storage function $V$ and finding the linear convergence rate $\alpha$.
\end{enumerate}

We will apply this methodology to analyze the DGDA update \eqref{eq: Dissipative}. Let $ z = [x,y]^T, \hat{z} = [\hat{x},\hat{y}]^T$ and $F(z_k) = (\nabla_x f(x_{k},y_{k}); -\nabla_y f(x_{k},y_{k}))$, and rewrite \eqref{eq: Dissipative} as in the form of \eqref{eq: dis dynamical system}:
\begin{equation}
    \begin{aligned}
        \begin{bmatrix}
    z_{k+1} \\
    \hat{z}_{k+1}
    \end{bmatrix}&= 
    \begin{bmatrix}
    z_k - \eta F(z_k)- \rho(z_k-\hat{z}_k) \\
    \hat{z}_{k} - \rho(\hat{z}_k -z_k)
    \end{bmatrix} \\
    &= 
    \begin{bmatrix}
        1-\rho & \rho \\
        \rho & 1-\rho 
    \end{bmatrix}\begin{bmatrix}
        z_k \\
        \hat{z}_{k}
    \end{bmatrix} + \begin{bmatrix}
        -\eta \\
        0
    \end{bmatrix}w_k 
    \end{aligned}
\end{equation}
where $w_k = F(z_k)$.

According to the previous discussion, the first step would be to choose a proper quadratic supply rate function $S$ such that $S( {\xi} _{k},w_k ) \leq 0 ,\forall k$, where ${\xi} _{k} =(z_k ;\hat{z}_{k})  $ that depends on the specific nonlinear term $w_k= F(z_k)$. According to the equations (7) in the work by Hu et al. (2017) \cite{hu2017dissipativity} and Lemma 6 from the research by Lessard et al. (2016) \cite{lessard2016analysis}, the following applies to the nonlinear operator $F(z_k)$ that meets the conditions specified in Assumption \ref{ass: Str}:
\begin{align}
    S(z_k,w_k) = \begin{bmatrix}
         z_k \\
        w_k
    \end{bmatrix}^T\begin{bmatrix}
        2\mu LI  &(-\mu+L)I \\
        (-\mu+L)I  & 2I \end{bmatrix}
        \begin{bmatrix}
         z_k \\
        w_k
    \end{bmatrix} \leq 0 
\end{align}
The conditions in Assumption \ref{ass: Str} are also commonly referred to as being L-smooth and m-strongly monotone, as can be found in related literature on variational inequality problems, such as the works by \cite{zhang2021unified,gidel2018variational}]. Therefore, we could easily extend the above LMI into the following supply rate function for DGDA updates, by augmenting the states ${\xi} _{k} =(z_k ;\hat{z}_{k})$:
\begin{equation}
    \begin{aligned}
      &S(\xi_k,w_k) = \\
    & \begin{bmatrix}
         z_k \\
        \hat{z}_{k} \\
        w_k
    \end{bmatrix}^T\begin{bmatrix}
        2\mu LI & 0 &(-\mu+L)I \\
        0 & 0&0\\
        (-\mu+L)I & 0 & 2I \end{bmatrix}
        \begin{bmatrix}
         z_k \\
        \hat{z}_{k} \\
        w_k
    \end{bmatrix} \leq 0   
    \end{aligned}
\end{equation} as a proper quadratic supply rate function $ S({\xi} _{k},w_k ) \leq 0$, whenever $w_k=F(z_k)$.

Finally, according to Theorem \ref{thm: Dissipative} and the above discussion, proving linear convergence reduces to finding a positive definite a matrix $P \in \mathbf{R}^{2(n+m) \times 2(n+m)}$, $\alpha \in [0,1)$ such that \eqref{eq: Dis LMI} is satisfied, where the problem parameters are given by
\begin{equation}
    \begin{aligned}
        &A =  \begin{bmatrix}
        1-\rho & \rho \\
        \rho & 1-\rho 
    \end{bmatrix}\otimes I, B = \begin{bmatrix}
        -\eta \\
        0
    \end{bmatrix}\otimes I, \\
    & X = \begin{bmatrix}
        2\mu L & 0 &(-\mu+L) \\
        0 & 0&0\\
        (-\mu+L) & 0 & 2 \end{bmatrix}\otimes I,
    \end{aligned}
\end{equation}

where $\otimes$ is the Kronecker product.
Due to the Kronecker structure of this problem, this is equivalent to solving an LMI problem of dimension 3 by 3, with design parameters $P=\bar P\otimes I$, with  $\bar P \in \mathbf{R}^{2 \times 2} $, $\alpha^2 \in [0,1)$, $\rho$ and $\eta$. 
Because this Linear Matrix Inequality is simple (3 by 3), it can be solved using analytical methods. This, in turn, results in a feasible solution for the LMI, denoted as follows:
\begin{equation}
    \begin{aligned}
        &\rho = \frac{1}{2} ,\qquad \eta = \frac{1}{L+\mu}, \\
        &\alpha^2 = \frac{3L^2 + 2L \mu + 3\mu^2 + \sqrt{(L+\mu)^4 + 16L^2\mu^2}}{4(L+\mu)^2}, \\
 & P = \begin{bmatrix}
      (L+\mu)^2 & 0\\
      0 & (L+\mu)^2
  \end{bmatrix}\otimes I.
    \end{aligned}
\end{equation}

After substituting the definition for condition number $\kappa:=L /\mu$, the convergence rate $\alpha^2$ simplifies to: 
\begin{align}
    \alpha^2 = 1 - \kappa^{-1} + \mathcal{O}\bigl((\frac{\mu}{L})^2 \bigr)
\end{align}

\subsection{Proof of Corollary \ref{cor:Str montonte comparison with known rates}} \label{Appendix: Str montonte comparison with known rates}
According to Theorem \ref{thm: Strongly monotone linear convergence}, the linear convergence rate estimate of DGDA is 
\begin{align}
     &\frac{3L^2 + 2L \mu + 3\mu^2 + \sqrt{(L+\mu)^4 + 16L^2\mu^2}}{4(L+\mu)^2} \\
     & = \frac{3\kappa^2+2\kappa+3+\sqrt{(\kappa+1)^4 + 16\kappa^2}}{4(\kappa+1)^2},
\end{align}
where $\kappa  = L/\mu$.

According to Theorem $6 \& 7$ \cite{azizian2020tight} and Theorem $4\&7$ \cite{mokhtari2020unified}, the standard known linear convergence rate estimate of EG and OGDA is 
\begin{align}
    1 - \frac{\mu}{4L} = 1- \frac{1}{4\kappa}.
\end{align}
By simple algebraic calculation, it can be shown that as a function of $\kappa$, when $\kappa \geq 2$, the following polynomial is always nonnegative, i.e.,
\begin{align}
     1- \frac{1}{4\kappa} - \frac{3\kappa^2+2\kappa+3+\sqrt{(\kappa+1)^4 + 16\kappa^2}}{4(\kappa+1)^2} \geq 0.
\end{align}




\end{document}